\documentclass{article}

\usepackage[utf8]{inputenc}
\usepackage{fullpage}
\usepackage{amssymb}
\usepackage{amsmath}
\usepackage{graphicx}
\usepackage{wrapfig}
\usepackage{cite}
\usepackage{subfigure}
\usepackage{makeidx}
\usepackage{multicol}
\usepackage{mathptmx}
\usepackage{mathrsfs}
\usepackage{mathtools}
\usepackage{amsthm}

\title{\textbf{A Control Framework for Socially-Optimal Emerging Mobility Systems}}
\author{Andreas A. Malikopoulos (amaliko@cornell.edu)\\ Professor\\ Cornell University}
\date{}

\begin{document}

\maketitle

\begin{center}
	\textbf{Abstract}
\end{center}
Connected and automated vehicles (CAVs) provide the most intriguing opportunity for enabling users to significantly improve safety and transportation efficiency by monitoring network conditions and making better operating decisions. CAVs, however, could alter the tendency to travel, eventually leading to a high traffic demand and causing rebound effects (e.g., increasing vehicle miles traveled). This chapter provides a control framework to distribute travel demand in a given transportation network, resulting in a socially optimal mobility system that travelers would be willing to accept. A “socially optimal mobility system” implies a mobility system that (1) is efficient (in terms of energy consumption and travel time), (2) mitigates rebound effects, and (3) ensures equity in transportation.

\section{Introduction}
\begin{figure}
	\centering
	\includegraphics[width=0.6\linewidth, keepaspectratio]{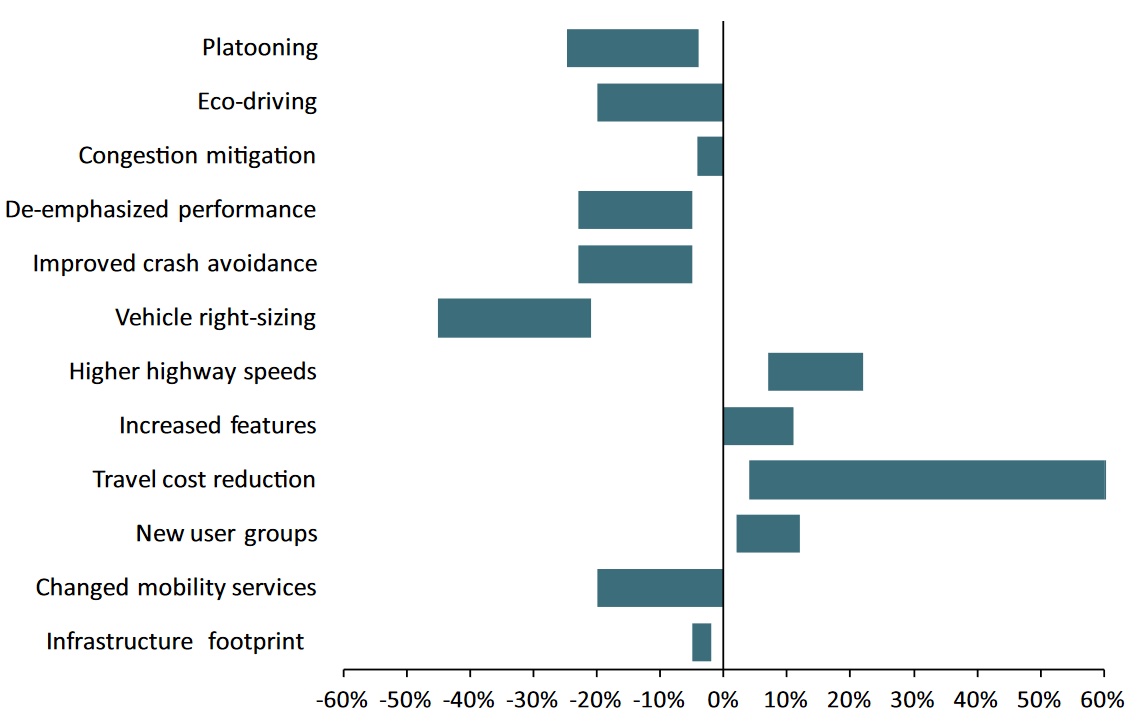} 
	\caption{Summary of estimated ranges of operational energy impacts of vehicle automation \cite{wadud2016}.}%
	\label{fig:impact}%
\end{figure}

Emerging mobility systems, such as connected and automated vehicles (CAVs) and shared mobility, provide significant opportunities to improve safety and reduce pollution, energy consumption, and travel delays \cite{zhao2019enhanced,Connor2020ImpactConnectivity}. CAVs are typical cyber–physical systems where the cyber component (e.g., data and shared information through vehicle-to-vehicle and vehicle-to-infrastructure communication) can aim at optimally controlling the physical entities (e.g., CAVs, and non-CAVs).
The {cyber-physical} nature of such emerging mobility systems requires large quantities of shared information through vehicle-to-vehicle and vehicle-to-infrastructure communication, which, in turn, are associated with significant technical challenges. This has given rise to a new level of complexity \cite{Malikopoulos2016} in modeling and control \cite{Ferrara2018}. 
It is expected that CAVs will gradually penetrate the market and interact with human-driven vehicles in ways that will improve safety and transportation efficiency over the next several years \cite{Malikopoulos2018d, Zhao2018CTA, cassandras2019b}. 
However, different levels of vehicle automation in the transportation network can have differing impacts on different transportation efficiency metrics \cite{wadud2016} ranging from an improvement of 45\% to a deterioration of 60\% as compared to the baseline of no automation (Fig. \ref{fig:impact}). 
Moreover, we anticipate that efficient  transportation and travel cost reduction might alter human travel behavior, causing {rebound effects.} For example, one impact of improving efficiency is a decrease in travel cost and a subsequent increase in the willingness to travel, as indicated in Fig. \ref{fig:impact}. The latter would increase overall vehicle miles traveled, which, in turn, might negate the benefits in terms of energy and travel time.

While several studies have shown the benefits of emerging mobility systems to reduce energy and alleviate traffic congestion in specific transportation scenarios, one {key question} that still remains unanswered is  ``How can we develop a mobility system that can enhance accessibility, safety, and equity in transportation without causing rebound effects, while also gaining the travelers' acceptance?"
This chapter aims to address this question using a control framework with two attributes:
\begin{enumerate}
	\item Distribute travel demand in a given transportation network, resulting in a {socially-optimal mobility system}.
	\item All travelers will willingly accept this distribution of travel demand. 
\end{enumerate}

In our exposition, the notion of a ``socially optimal mobility system'' implies that the mobility system is efficient (in terms of energy consumption and travel time), mitigates rebound effects, and ensures {equity} in transportation.
Specifically, the control framework involves the following two research steps: 

\textbf{Step 1:} Aggregate the preferences of the travelers into a collective, system-wide set of recommendations, e.g., routing choices, modes of transportation, while the private information of the travelers is not publicly known; and 

\textbf{Step 2:} Develop control technologies allowing CAVs to navigate automatically and to co-exist with human-driven vehicles (HDVs) safely and efficiently in a mixed traffic environment.

Step 1 will identify the new congestion patterns of an optimized mobility system. In contrast, Step 2 will improve transportation efficiency under safety guarantees in response  to the new levels of imposed travel demand. 
More specifically, the two research steps address the following challenges: 
In Step 1, we  formulate and solve an optimization problem in which the decision variables are the optimal routes and the selection of a transportation mode for all travelers so as to maximize a social utility function. A suitable approach to address this problem is considering a decentralized traveler decomposition using mechanism design theory \cite{hurwicz2006}. In mechanism design, we are concerned with implementing system-wide optimal solutions to problems involving multiple agents -- in this case, travelers -- each with private information about preferences, e.g., individual tolerance to traffic delay, the value of time and money, and preferred travel time.  In this context, a social planner faces the problem of aggregating the travelers' preferences into a collective, system-wide decision when the travelers' private information is not publicly known. Thus, mechanism design entails the social planner solving an optimization problem with {incomplete information.}
\\

In Step 2,  we need to have an optimal coordination framework which, in conjunction with large amounts of data from vehicles and the infrastructure, will improve safety and efficiency  in a mixed traffic environment consisting of CAVs and human-driven vehicles. In particular, the coordination framework will be used by CAVs to navigate optimally (in the sense of minimizing travel times and energy consumption) while also guaranteeing safety in different traffic settings, e.g., crossing a signal-free intersection without stopping, merging at roadways or a roundabout, and executing automated passing maneuvers. This step will ensure that CAVs co-exist and interact {safely} with other human-driven vehicles and pedestrians.
\\

The one feature that {sharply distinguishes} the proposed framework from previous approaches reported in the literature to date is that  it considers simultaneous optimization of the {travel demand} (Step 1) and {transportation network efficiency} (Step 2).

\subsection{A Brief Review of Mechanism Design}

Mechanism design theory is concerned with implementing system-wide optimal solutions to problems involving multiple agents, each with private information about preferences. It can be viewed as the art of designing the rules of a game to achieve a specific desired outcome. A popular example of mechanism design, taken from [2], is a cake-cutting problem involving a mother with two kids. The mother must design a mechanism to make her kids share a cake equally. The centralized solution of this problem is for the mother to slice the cake equally and give a slice to each kid. The decentralized mechanism design solution is to (1) allow one kid to slice the cake into two pieces and (2) allow the other kid to determine who gets which piece. This mechanism achieves the desired outcome of the kids sharing the cake equally without the mother’s intervention. In this example, the mother is called the “social planner.” In the context of mechanism design, the social planner typically faces the problem of aggregating the preferences of multiple agents into a collective, system-wide decision when the agents' private information might not be publicly known. Thus, mechanism design entails the social planner to solve an optimization problem with incomplete information. The underlying structure used is to induce a game among the agents so that the desired system-wide solution is implemented in an equilibrium of the induced game. 
\\

Numerous economists and mathematicians have rigorously studied the mechanism design theory since the 1950s to provide insights and solutions to different economic topics. The theory started with the seminal contributions of Leonid Hurwicz and Jacob Marschak, who were interested in resource allocation problems and controlling (through incentives) individual agents. In parallel, Kenneth Arrow, G\'{e}rard Debreu, and Herbert A. Simon also worked on problems with incomplete information and how to bound rationality. In 1961, William Vickrey's seminal work on auctions was published, paving the way for Hurwicz's theoretical framework to be applied to design incentives based on the agents' information for a simple yet formidable problem of an auction. Arrow and Debreu's work was instrumental in establishing the interconnected relation of information and decision-making and its role in influencing behavior in the efficient allocation of limited resources. Much later in the 1970s and 1980s, Peter Diamond, Oliver Hart, Jean-Jacques Laffont, Eric Maskin, James Mirrlees, and Sherwin Rosen worked independently on principal-agent problems focusing on how one can design a contract between a principal (e.g., institution corporation) and a rational agent efficiently. As a continuation of Vickrey's work, Ronald Coase, Jerry R. Green, Theodore Groves, and John Ledyard made significant contributions to the design of incentives for public goods problems. Furthermore, Roger Myerson, Paul Milgrom, and Robert Wilson expanded the Vickrey auction to address more complicated scenarios and complex problems.
The theory of mechanism design represents the confluence of microeconomics \cite{Laffont1996,Clarke1971,Groves1975,Maskin1999,Green1977,Groves1973,Groves1977,Maskin2008,Maskin2000,Reichelstein1988,Armstrong1996,Mathevet2010}, and social choice theory \cite{acemoglu2008political,acemoglu2010dynamic}, while it equally draws from auction \cite{Vickrey1961}, optimization \cite{bandi2012tractable,Bandi2014}, and game theory \cite{Roughgarden2019,Dasgupta1979,Maskin2002,Richter2019,Condorelli2012,Thomson1996}.

\section{Outline of the Control Framework}
\label{novel}

It is expected that CAVs will gradually penetrate the market and disturb travelers' behaviors along with their mobility tendencies, resulting in unintended consequences (rebound effects). These consequences may be additional energy use and greenhouse gas emissions, challenging equity in transportation,  and significant alterations in the density of urban areas \cite{anderson2014,fagnant2015,litman2020}. Consequently, we are compelled to reassess the relationship between mobility and social life \cite{sheller2000,bissell2020} and provide solutions that consider the {decision-making} processes of travelers. 

\begin{figure}
	\centering
	\includegraphics[width=0.8\linewidth, keepaspectratio]{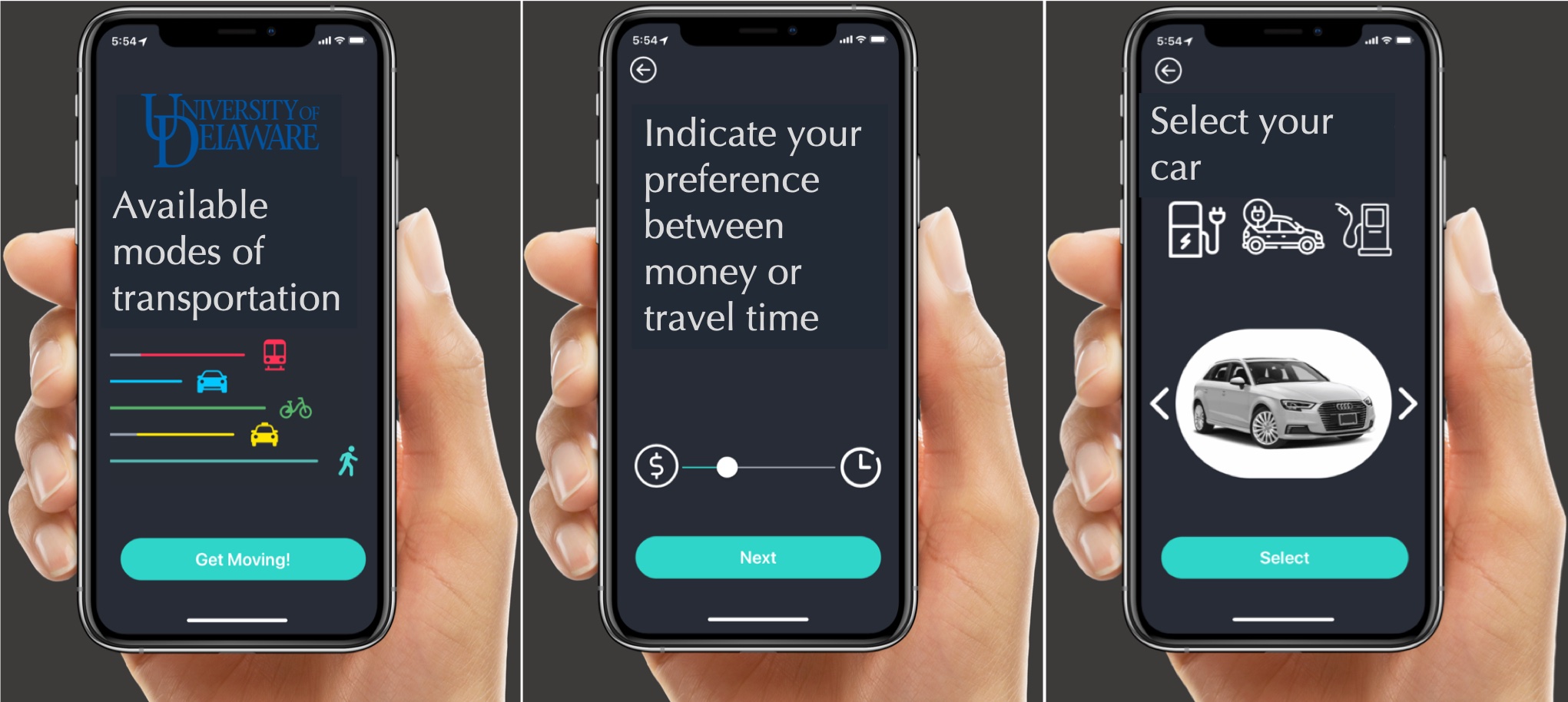} 
	\caption{The app of the proposed framework.}%
	\label{fig:app}%
\end{figure}

Our aim is to develop a holistic and rigorous framework to capture the {societal} impact of emerging mobility systems \cite{chremos2021MobilityGame} and provide solutions that mitigate any potential rebound effects, e.g., increased vehicle miles traveled, increased travel demand, and empty vehicle trips, while enhancing accessibility, safety, and equity in transportation.
In our approach, we consider a finite set of travelers who seek to travel in a given transportation network of a city (or a big metropolitan area) where a central authority (the {social planner} in our exposition) seeks to ensure the efficient allocation and operation of the different modes of transportation available in the transportation network of the city. We call these different modes {``mobility services."} A few examples are privately-owned CAVs and human-driven vehicles, shared mobility vehicles (e.g., Uber, Lyft), bicycles, and public transit (e.g., train, bus, light rail). The travelers make a request via a smartphone app (Fig. \ref{fig:app})  to use a service to satisfy their mobility needs, i.e., desired origin destination. The social planner, which  can be visualized as a central computer, compiles all travelers' origin-destination requests and other information, e.g., preferred travel time, value of time, to provide a travel recommendation to each traveler with the goal of achieving a {socially-optimal} mobility system.
In this context, the social planner formulates and solves an optimization problem in which the decision variables are the mobility services that maximize the travelers' social welfare (utilities).
A suitable approach to address this problem is to consider a decentralized decomposition of the travelers using mechanism design theory \cite{borgers2015,hurwicz2006}. In mechanism design, we are concerned with how to implement system-wide optimal solutions to problems involving multiple agents -- in this case, travelers -- each with private information. Thus, mechanism design entails the social planner solving an optimization problem with incomplete information.
\\

The social planner's collective recommendations must achieve three objectives: {(1)} respect and satisfy the travelers' preferences, {(2)} guarantee equity in travel recommendations, and {(3)} ensure that no mobility system will become congested. We assume the city supports connected and automated mobility technologies on public transit infrastructure and a transportation network. Consequently, the social planner is {fully aware} of the system's capabilities and the network's capacity. In other words, the social planner is fully capable of computing the maximum capacity of each mobility service and the associated costs aimed at providing travel recommendations to all travelers.
The social planner's objective is to {design} appropriate monetary incentives, e.g., tolls, fares, subsidies, to guarantee the realization of the desired outcome, i.e., maximize the social welfare of all travelers. The travelers will, in turn, accept or reject the social planner's recommendations.

\section{Socially-Optimal Management of Travel Demand}
\label{task:1.1}
As we move to increasingly complex emerging mobility systems with an expanded feature space, fundamentally new approaches are needed to understand the impact on system behavior \cite{Malikopoulos2016b}.  The approaches reported in the literature to date have considered emerging mobility systems without deliberating on {human decision-making} and {perception.} 
To develop and operate a socially-optimal mobility system, {technological} and {information management} innovations need to be integrated with the {social} dimensions to ensure adoption by the drivers, travelers, and the public. 
\\

The standard approach to alleviate congestion in transportation has been the management of travel demand. Some approaches have considered congestion pricing/tolling \cite{vickrey1969} while others have considered the application of mechanism design to provide a solution to individual route selection under different congestion traffic scenarios. The theory of mechanism design was developed for the implementation of system-wide optimal solutions to problems involving multiple rational agents, each with private information and conflicting interests \cite{mas_colell1995}. It can be viewed as the art of designing the rules of a game to achieve a specific desired outcome. Mechanism design has broad applications spanning different fields, including microeconomics, social choice theory, and control engineering. 
Applications in engineering include communication networks \cite{renou2012}, social networks \cite{Dave2020SocialMedia}, transportation routing \cite{bian2019}, online advertising \cite{kakade2013}, smart grid \cite{samadi2012}, multi-agent systems \cite{shoham2008}, and resource allocation problems \cite{zou2015}.
Due to the economic nature of congestion in transportation, auctioning has also been proposed \cite{iwanowski2003} to create a market of tolls in a network of roads. Auctions are processes for allocating goods among bidders, so the challenge of auction design can only be understood by studying the demands of the participants \cite{mishra2009}. Auction design has been the focus of significant results on multi-object auctions and matching market problems \cite{demange1986,gale1962}. On the one hand, auctions have been proposed to design pricing schemes with tolls in a network of roads, leading to a spark of studies in auctioning techniques \cite{iwanowski2003,teodorovic2008,vasirani2011,raphael2015,olarte2018}. On the other hand, this approach has significant limitations: (1) the implementability of auction-based tolling on highways is not straightforward due to the dynamic and fast-changing nature of transportation systems; (2) it is also uncertain how the public, e.g., travelers, passengers, drivers, will respond concerning toll roads in an auction setting. Understanding the travelers' interests, i.e., willingness-to-pay, the value of time,  and the impacts on different sociodemographic groups, becomes imperative for a socially efficient design of an emerging mobility system.
\\

In our approach, we seek to ``design" a socially optimal mobility system that assigns mobility services to a finite group of travelers by considering their personal travel preferences.  By a socially optimal mobility system, we mean a mobility system that is {(1)} efficient (in terms of energy consumption and travel time), {(2)} mitigates rebound effects, and {(3)} ensures equity in transportation. We want to ensure that the emerging mobility system is {incentive compatible} (travelers always report their personal travel preferences truthfully), {individually rational} (travelers always benefit from voluntarily participating in the system), and {weakly budget balanced} (the system always generates revenue from each traveler that can be used for maintenance of the road network, etc). 
\\

Our mobility system is managed by a {social planner} who aims to allocate $m \in \mathbb{N}$ mobility services to $n \in \mathbb{N}$ travelers, where $n \geq m$. We denote the nonempty set of travelers by $\mathcal{I} = \{1, 2, \dots, n\}$ and the nonempty set of mobility services by $\mathcal{J} = \{1, 2, \dots, m\}$. For example, each mobility service $j \in \mathcal{J}$ can either represent privately-owned CAVs and HDVs, shared mobility vehicles (e.g., Uber, Lyft), bicycles, and public transit (e.g., train, bus, light rail). For our purposes, we can think of $\mathcal{G} = (\mathcal{V}, \mathcal{E})$ representing a smart city network with a road and public transit infrastructure. 
Next, a traveler $i \in \mathcal{I}$ seeks to travel using these mobility services in a transportation network from their current location $o_i \in \mathcal{V}$ to their desired destination $d_i \in \mathcal{V}$. So, each traveler $i \in \mathcal{I}$ is associated with the origin-destination pair $(o_i, d_i)$ represented by an undirected multigraph $\mathcal{G} = (\mathcal{V}, \mathcal{E})$, where each node in $\mathcal{V}$ represents a different neighborhood, and each link $e \in \mathcal{E}$ represents a sequence of city roads and a public transit connection. 
On the other hand, each mobility service $j \in \mathcal{J}$ is associated with several different links. Therefore, we do not have to limit the number of mobility services that connect any origin $o_i$ to any destination $d_i$. Furthermore, 
Each traveler $i \in \mathcal{I}$ seeks to travel in the network $\mathcal{G}$ with only one mobility service $j \in \mathcal{J}$ that satisfies their origin-destination pair $(o_i, d_i)$ while each service $j \in \mathcal{J}$ can be used by multiple travelers.
\\

The traveler-service {assignment} is a vector $\mathbf{a} = (a_{1 1}, \dots, a_{i j}, \dots, a_{n m}) = (a_{i j})_{{i \in \mathcal{I}}, j \in \mathcal{J}}$, where $a_{i j}$ is a binary variable of the form: $a_{i j} =1$, if $i \in \mathcal{I}$ is assigned to $j \in \mathcal{J}$, and $a_{i j} =0$, otherwise.
The set of travelers with the same origin-destination pair is 
\begin{align}
	\mathcal{I}_k = \{i \in \mathcal{I} \; | \; (o_i, d_i) = (o_k, d_k)\}, ~k = 1, 2, \dots, N, 
\end{align}
where $N \in \mathbb{N}$ is the number of sub-classes over the complete set of travelers, i.e., $\mathcal{I} = \bigcup_{k = 1} ^ N \mathcal{I}_k$.
The justification of the introduction of  $\mathcal{I}_k $ is that in an emerging mobility system, we can acquire verifiable location data of travelers either by using a global positioning system or estimating the average number of travelers using public transit \cite{coleri2004}.
By partitioning the set of travelers in $N \in \mathbb{N}$ sub-classes, the traveler-service assignment of sub-class $\mathcal{I}_k$ is given by $\mathbf{a}_k = (a_{i j})_{{i \in \mathcal{I}_k}, j \in \mathcal{J}}$.
Since each mobility service $j \in \mathcal{J}$ can be utilized by multiple travelers simultaneously, we introduce a metric of ``co-travelers." Hence, for any traveler, $i \in \mathcal{I}_k$, $k = 1, \dots, N$, the number of co-travelers having the same origin-destination pair and using the same service $j \in \mathcal{J}$ is computed by	
\begin{align}
	\psi_i(\mathbf{a}_k) = \sum_{\ell \in \mathcal{I}_k} a_{\ell j} - a_{i j}.
\end{align}

A traveler $i \in \mathcal{I}_k$, $k = 1, \dots, N$, is characterized by a tuple of {personal travel preferences}, denoted by $\pi_i$, and given by
$\pi_i = (\theta_i, \eta_{i j}, \delta_{i j}),$
where $\theta_i \in \mathbb{R}_{\geq 0}$ is the preferred travel time, $\eta_{i j} \in \mathbb{N}$ is the maximum preferred number of co-travelers with the mobility service $j$, and $\delta_{i j} \in [0, 1]$ is a monetary value that traveler $i$ is willing to pay (or accept as a compensation) when using mobility service $j$.
Naturally, a traveler $i$'s preferred travel time $\theta_i$ is a non-negative real number and represents how fast traveler $i \in \mathcal{I}_k$ wishes to reach their destination. Similarly, traveler $i$'s preferred number of co-travelers $\eta_{i j}\in\mathbb{N}$ represents the maximum tolerable number of other travelers using mobility service $j$, and $\delta_{i j}$ represents what potentially drives a traveler's behavior, i.e., the value of time. In particular, $\delta_{i j}$ represents a monetary value that traveler $i \in \mathcal{I}_k$ is willing to pay to save time (or accept as compensation) for service $j \in \mathcal{J}$.
For each traveler $i \in \mathcal{I}_k$, $\pi_i $ is considered {private information,} known only to traveler $i \in \mathcal{I}_k$. Hence, $(\pi_i)_{i \in \mathcal{I}_k}$, $k = 1, \dots, N$ are unknown information to the social planner. This is one of the {key challenges} in the proposed mobility system, i.e., ``how do we {incentivize} the travelers to be {truthful} and {elicit} the private information needed to provide a socially optimal solution for the whole system?" The answer to this question is given in Section \ref{sec2}.
\\

Next, we introduce an {``inconvenience"} metric for any traveler $i \in \mathcal{I}_k$ using any mobility service $j \in \mathcal{J}$. Quantitatively, the inconvenience metric can represent the extra monetary value of travel dis-utility from any costs, travel delays, or violation of personal preferences caused by using the mobility service $j \in \mathcal{J}$. The {mobility inconvenience metric} for traveler $i \in \mathcal{I}_k$, $k = 1, \dots, N,$ assigned to service $j \in \mathcal{J}$ is a continuous function 
\begin{align}
	\phi_i ( \pi_i, \tilde{\theta}_i(\mathbf{a}_k), \psi_i(\mathbf{a}_k) ) \in \mathbb{R}_{\geq 0},
\end{align}
where $\pi_i$ is the tuple of the personal travel preferences of $i \in \mathcal{I}_k$, $\tilde{\theta}_i(\mathbf{a}_k) \in \mathbb{R}_{\geq 0}$ is the experienced travel time, and $\psi_i(\mathbf{a}_k)$ is the number of co-travelers.
Note that the mobility inconvenience metric $\phi_i$ strictly increases when $\tilde{\theta}_i(\mathbf{a}_k)$ and/or $\psi_i(\mathbf{a}_k)$ increases. That is because, from a modeling perspective, traveling with time delays or during peak times can cause significant inconveniences to any traveler $i \in \mathcal{I}_k$.
\\

Next, a traveler's satisfaction is captured by a valuation function $v_i(\mathbf{a}_k)$, which can reflect the traveler's willingness-to-pay for their travel, i.e.,
\begin{equation}\label{EQN:mobility-valuation}
	v_i(\mathbf{a}_k) = \bar{v}_i - \phi_i \left( \pi_i, \tilde{\theta}_i(\mathbf{a}_k), \psi_i(\mathbf{a}_k) \right),
\end{equation}
where $\bar{v}_i \in \mathbb{R}_{> 0}$ is the value gained by traveler $i \in \mathcal{I}_k$ when their origin-destination pair $(o_i, d_i)$ is satisfied using mobility service $j \in \mathcal{J}$ without any travel delays or travel inconveniences, i.e., $\tilde{\theta}_i(\mathbf{a}_k) = \theta_i$ and $\psi_i(\mathbf{a}_k) = 0$. We call $\bar{v}_i$ the \emph{maximum willingness-to-pay}. Naturally, for any traveler $i \in \mathcal{I}_k$ and any service $j \in \mathcal{J}$, we have $v_i(\mathbf{a}_k) \in [0, \bar{v}_i]$, where $v_i(\mathbf{a}_k) = 0$ means that traveler $i \in \mathcal{I}_k$ is unwilling to utilize mobility service $j \in \mathcal{J}$. Although our analysis will treat the satisfaction function $v_i(\mathbf{a}_k)$ in its most general form, given by \eqref{EQN:mobility-valuation}, one can explicitly define it as follows: 
\begin{equation}
	v_i(a_{i j}) =
	\begin{cases}
		\bar{v}_i, & \text{if } \phi_i = 0, \\
		\lambda_i \cdot \bar{v}_i, & \text{if } \phi_i = (1 - \lambda_i) \cdot \bar{v}_i, \\
		0, & \text{if } \phi_i = \bar{v}_i,
	\end{cases}
\end{equation}
where $\lambda_i \in (0, 1)$ is a discount rate. When $\phi_i = 0$, we say that traveler $i$'s personal travel requirements $\pi_i$ are satisfied using the CAV mobility service without any co-travelers. When $\phi_i = \bar{v}_i$, then we say that traveler $i$'s personal travel requirements are not satisfied. Lastly, when $\phi_i = (1 - \lambda_i) \cdot \bar{v}_i$, we say that traveler $i$'s personal travel requirements $\pi_i$ are satisfied using any mobility service.

The {total utility} $u_i(\mathbf{a}_k)$ of traveller $i \in \mathcal{I}_k$, $k = 1, \dots, N,$ is given by
\begin{align}
	u_i(\mathbf{a}_k) = v_i(\mathbf{a}_k) - p_i(\mathbf{a}_k),
\end{align}
where $v_i(\mathbf{a}_k)$ is the maximum willingness-to-pay and $p_i(\mathbf{a}_k)$ is the {mobility payment} traveler $i \in \mathcal{I}_k$ is required to make to use service $j \in \mathcal{J}$ (e.g., pay road tolls or buy a public transit ticket).
The {operating cost} of service $j \in \mathcal{J}$, denoted by $r_j \in \mathbb{R}_{> 0}$, is given by
\begin{align}
	r_j(\mathbf{a}_k) = \sum_{i \in \mathcal{I}_k} r_{i j}(a_{i j}),
\end{align}
where $r_{i j}(a_{i j}) \in \mathbb{R}_{> 0}$ is traveler $i$'s corresponding share of the operating cost of vehicle $j \in \mathcal{J}$.
Intuitively, the operating cost $r_{i j}$ captures traveler $i$'s fair share of the costs of mobility service $j \in \mathcal{J}$. These costs can be associated with fuel/energy consumption, drivers' labor reimbursement, and environmental impact. Moreover, the operating cost $r_{i j}$ can be thought as the minimum acceptable payment, i.e., $p_i \geq r_{i j}$, by traveler $i \in \mathcal{I}_k$ using service $j \in \mathcal{J}$.
\\

In the modeling framework described above, we impose the following {assumptions: 
	
	\textbf{Assumption 1:}	For all sub-classes $\mathcal{I}_k$, $k = 1, \dots, N$, $N \in \mathbb{R}$, any traveler $i \in \mathcal{I}_k$ is modeled as a selfish decision-maker with private information $\pi_i = (\theta_i, \eta_{i j}, \delta_{i j})$. Traveler $i$'s objective is to maximize their total utility 
	\begin{align}
		u_i(\mathbf{a}_k) = v_i(\mathbf{a}_k) - p_i(\mathbf{a}_k).
	\end{align}
	
	\textbf{Assumption 2:} For any sub-class $\mathcal{I}_k$, $k = 1, \dots, N$, and for any traveler $i \in \mathcal{I}_k$ the maximum willingness-to-pay realized from using any mobility service $j \in \mathcal{J}$ must outweigh the operating costs, i.e., for all $i \in \mathcal{I}_k$, we have $\bar{v}_i > r_{i j}(a_{i j})$, where $r_{i j}(a_{i j})$ is the specific operating cost imposed by traveler $i \in \mathcal{I}_k$ using mobility service $j \in \mathcal{J}$.
	
	The first assumption essentially indicates that each traveler is selfish in the sense that they are only interested in their own well-being. In economics, such behavior is called ``strategic" since agents attempt to misreport or lie about their private information to the social planner if that means higher individual benefits. The second assumption implies it is always beneficial to travel when it can be guaranteed that there will be no travel delays or inconveniences.
	\\
	
	In the proposed framework, travelers request (via a smartphone app) in advance a travel recommendation from the social planner that satisfies their origin-destination. Given the travelers' origin-destination pairs, the social planner distributes all travelers to different sub-classes. Thus, travelers from the same neighborhood have the same origin. Similarly, travelers going to the same neighborhood have the same destination. The social planner's task is to satisfy all travel requests and provide recommendations to the travelers, e.g., which mobility service to use. Hence, we are interested in minimizing the travel inconvenience of all travelers and the operating costs, which is equivalent to maximizing the utility of each traveler. 
	Thus, the social planner formulates the following  {optimization problem} for each sub-class $\mathcal{I}_k$, $k = 1, \dots, N,$ 
	\begin{align}
		J_1(\mathbf{a}_k) = \min_{\mathbf{a}_k}[\omega_1 \sum_{i \in \mathcal{I}_k} \phi_i ( \pi_i, \tilde{\theta}_i(\mathbf{a}_k), \psi_i(\mathbf{a}_k) ) + \omega_2\sum_{j \in \mathcal{J}} r_j(\mathbf{a}_k)],\label{problem1}
	\end{align}
	where $\omega_1$ and $\omega_2$ are factors that normalize the terms.
	The problem is subject to the following {constraints: (1)} each traveler $i \in \mathcal{I}_k$ is assigned at most one mobility service, {(2)} service $j$ must not exceed its maximum usage capacity, and {(3)} the traveler's assignment impose equity in transportation \cite{chremos2022Equity} using the mobility equity metric \cite{Bang2023mem}.

	\subsection{Desired Properties of the Mobility System}
	\label{sec2}
	
	The problem in \eqref{problem1} is a mixed-integer programming model, and standard algorithmic approaches exist to find its global optimal solutions or, in worst-case scenarios, their approximations. Note, though, that these approaches assume complete information of all parameters and variables in the model. Such an assumption is unreasonable to expect from strategic decision-makers. Thus, in our framework, travelers are not expected to report their private information truthfully. This turns our problem into a {preference elicitation problem}. Next, we discuss how we can elicit the necessary private information of all travelers using monetary incentives in the form of mobility payments (e.g., tolls, fares, fees).
	\\
	
	Since $(\pi_i)_{i \in \mathcal{I}_k}$, $k = 1, \dots, N,$ are unknown information for the social planner, the question becomes ``How do we {incentivize} the travelers to be {truthful} and {elicit} the private information needed to provide a socially-optimal solution for the whole system?" 
	We employ the mobility payments $p_i(\mathbf{a}_k)$ reported in \cite{chremos2020MobilityMarket} for each traveler $i \in \mathcal{I}_k$ to look  similar to the celebrated Vickrey-Clarke-Groves (VCG) mechanism \cite{Vickrey1961,Clarke1971,Groves1973}. Our  payments are significantly different within our context as we have considered the travelers' personal travel preferences and introduced the capacity constraints for each mobility service along with the equity in transportation constraint. In contrast, the VCG mechanism does not have any constraints, and thus, the proposed mobility system is considerably different from the VCG mechanism. 
	
	Therefore, a direction for future research should show that the mobility system satisfies (1) {incentive compatibility} and (2) {individual rationality} despite our departure from VCG.
	Incentive compatibility means that all travelers are incentivized to report their personal travel preferences truthfully regardless of what other travelers report. 
	Individual rationality implies that all travelers voluntarily participate in the  mobility system in the most potent form. Informally, we compare the utility of a traveler $i \in \mathcal{I}_k$ under two possible scenarios: traveler $i \in \mathcal{I}_k$ participates in the mobility market, and traveler $i \in \mathcal{I}_k$ rejects any travel recommendations from the social planner and simply uses their self-owning vehicle (CAV or conventional vehicle).
	The mobility system is individually rational if for any traveler $i \in \mathcal{I}_k$, we have	$u_i(\mathbf{a}_k) \geq u_i(\widehat{\mathbf{a}}_k),$
	where $\widehat{\mathbf{a}}_k$ denotes the traveler-assignment in which traveler $i \in \mathcal{I}_k$ rejects social planner's travel recommendations and instead uses their  self-owning vehicle. Finally, future research should show that the proposed mobility system is guaranteed to generate revenue from each traveler. This revenue can be used to maintain and update the infrastructure of the city's transportation network over the years.
	\\
	
	One particular limitation of the proposed framework is that we consider travelers' preferences to be static. 
	This implies that if any preferences of the travelers change, then the social planner would have to recompute the solution of the optimization problem in the mobility system to get an updated traveler-service assignment.

	\section{Optimal Coordination of Connected and Automated Vehicle with Human-Driven Vehicles}
	
	This section addresses the optimal coordination of CAVs with the neighbor human-driven vehicles (HDVs). We consider a team consisting of CAVs and HDVs that is about to encounter a given traffic scenario (e.g., crossing a signal-free intersection, merging at roadways or a roundabout, cruising in congested traffic, passing through a speed reduction zone, and lane-merging or passing maneuvers) with the {common objective} to {coordinate} in this scenario and avoid stop-and-go driving \cite{Zhao2018ITSC,Malikopoulos2020,chalaki2020TITS,mahbub2020Automatica-2,chalaki2020ICCA}. The implications of the latter are that the vehicles do not have to come to a full stop, thereby conserving momentum and fuel while also improving travel time. For example, consider a signal-free intersection (Fig. \ref{fig:3}) with a team of CAVs and HDVs. The region at the center of the intersection, called \textit{merging zone}, is the area of potential lateral collision of the vehicles. The intersection has a \textit{control zone} inside of which the CAVs can communicate with each other. The objective of the team of CAVs and HDVs is to cross the intersection without the use of traffic lights, without creating congestion, and under the hard safety constraint of collision avoidance. We should emphasize that the proposed  framework can be applied to any traffic scenario. We use an intersection here just as a reference for our exposition.
	
	\begin{figure}
	\centering
	\includegraphics[width=0.5\linewidth, keepaspectratio]{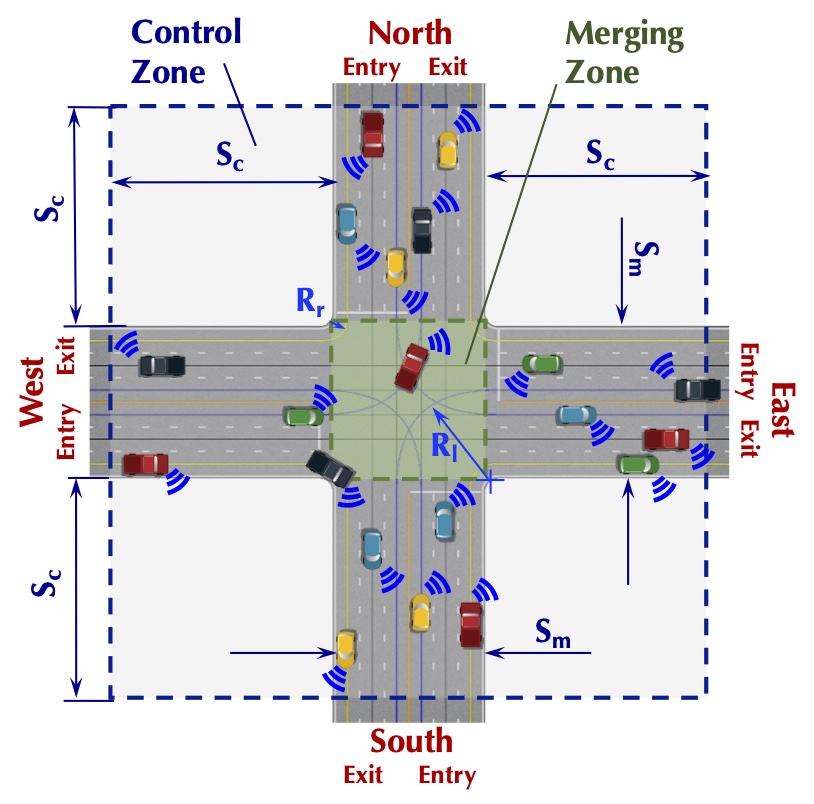} 
	\caption{
		A signal-free intersection with human-driven and connected automated vehicles.
	}%
	\label{fig:3}%
\end{figure}	
	
	We  model the communication between the team members with the {word-of-mouth communication} structure that we have previously developed \cite{Dave2020,Dave2020a}. In a word-of-mouth communication structure, every member of the team communicates with her neighbors with delays in communication. This is a {non-classical} information structure \cite{Malikopoulos2021} where the topological and temporal restrictions in communication mean that information propagates slowly through the team members. 
	\\
	
	In our modeling framework, we  consider a number of $K \in \mathbb{N}$ members in the team  of CAVs and HDVs with a {decentralized information} structure.
	At time $t = 0,1,\ldots,T$, $T\in\mathbb{N}$, the state of the team $X_t$ takes values in a finite set $\mathcal{X}_t$ and the decision $U_t^k$ associated with the team member $k \in \mathcal{K}$, $\mathcal{K}=\{1,\ldots,K\}$, takes values in a finite set $\mathcal{U}_t^k$. Let ${U}_t^{1:K}=(U_t^1,\ldots,U_t^K)$ be the team's decision at time $t$. Starting at the initial state $X_0$, the evolution of the team is described by the state equation
	$X_{t+1}=f_t\left(X_t,U_t^{1:K},W_t\right)$, where $W_t$ is a random variable corresponding to the external, uncontrollable disturbance to the team's mission that takes values in $\mathcal{W}$. The sequence $\{W_t: t=0,\ldots,T\}$ is a sequence of independent random variables which is also independent of the initial state $X_0$. 
	At time $t = 0,1,\ldots,T$, every team member $k\in\mathcal{K}$ makes an observation $Y_t^k$ according to the equation $Y_t^k=h_t^k(X_t,V_t^k)$, where $V_t^k$ is a random variable corresponding to the noise of the observation and takes values in the finite set $\mathcal{V}^k$. 
	\\
	
	To capture the delay typically encountered between vehicle-to-vehicle and vehicle-to-infrastructure communication, we consider that the team has $n$-step delayed information sharing. Namely, at time $t$, team member $k$ observes $Y_t^k$, and the $n$ steps past observations $Y_{t-n}^{1:K}$ and decisions $U_{t-n}^{1:K}$ of all team members. Thus, at time $t$ the information available at the team member $k$ is $(\Delta_t, \Lambda_t^k)$, where $\Delta_t\colon= (Y_{t-n}^{1:K}, U_{t-n-1}^{1:K})$, $Y_{t-n}^{1:K}=\{Y_{t-n}^{1},\ldots,Y_{t-n}^{K}\}$, $U_{t-n-1}^{1:K}=\{U_{t-n-1}^{1},\ldots,U_{t-n-1}^{K}\}$, is the information known to all team members and $\Lambda_t^k\colon= (Y_{t-n+1:t}^{k},$ $U_{t-n+1:t-1}^{k})$ is the additional information known at the team member $k\in\mathcal{K}$ only. 
	Note that the $n$-step delayed information sharing can  also be {asymmetric,} i.e., for each subsystem $k\in\mathcal{K}$, $Y_{t-n_i}^{k}$, $U_{t-n_i}^{k},$ where $n_i$, $i\in\mathcal{K}$, are constant but not necessarily the same for each  $k$.
	The collection $\{(\Delta_t, \Lambda_t^k);$ $ t=0,\ldots,T\}$, is the {information structure} of the team and captures who knows what about the team and when. Note that  $\Lambda_t^k$ of each robotic vehicle $k\in\mathcal{K}$ is {``richer"} than the information $\Lambda_t^j$ of each human-driven vehicle $j\in\mathcal{K,}$ since the observation $Y_{t}^{k}$ of $k$ also includes information from other CAVs, whereas HDVs can make only local observations, e.g., distance from the preceding vehicle, etc.
	\\
	
	Let $\mathcal{D}_t$ be the space of all possible realizations of $\Delta_t$, and $\mathcal{L}^k$ be the space of all possible realizations of  $\Lambda_t^k$. 
	The team member $k$ makes a decision according to a control law $g_t^k$, i.e., $U_{t}^{k}=g_t^k(\Delta_t, \Lambda_t^k)$. 
	The problem for each team member is to derive its {optimal control law} $g_t^k$ such that the collective {optimal strategy} ${g}=\{g_t^k; k\in\mathcal{K}; t=1,\ldots,T\}$ will enable the team to pass through the traffic scenario without stop-and-go driving. The latter can be modeled as the minimization of an expected total cost 
	\begin{align}\label{problem1a}
		J({g})=\mathbb{E}^g\big\{\sum_{t=1}^T c_t(X_t, U_t^{1:K})\}, 
	\end{align}
	where $c_t(X_t, U_t^{1:K})$ corresponds to travel delay of the team and the expectation is with respect to the joint probability of the random variables designated by the choice of ${g}$. 
	Deriving the solution of \eqref{problem1a} has the following conceptual difficulties:
	(1) the functional optimization problem of selecting a sequence of strategies is not trivial, as the set of the class of strategies is infinitely large and (2) the domain of the control strategies given by the information $\{(\Delta_t, \Lambda_t^{1:K}); t=0,\ldots,T\}$ increases with time, causing significant implications on storage requirements and real-time implementation.
	These difficulties can be circumvented by using the conditional probability of the state given the data available as a sufficient statistic, i.e., $\Pi_t\colon=\mathbb{P}(X_{t-n}|\Delta_t)$.
	The conditional probability is called {information state}, $\Pi_t$, and takes values in a time-invariant space.
	Using the information state can help us restrict our attention to control strategies in a {time-invariant domain.} 
	Such results, where data that are increasing with time are ``compressed" to a sufficient statistic taking values in a time-invariant space, are called {structural results.} 
	The structural results are related to the concept of separation, namely, the information state does not depend on the control strategy. This has been called a one-way separation between estimation and control. An essential consequence of this separation is that for any given choice of control strategies until time $t$ and a given realization of the system variables till time $t$, the information states at future times do not depend on the choice of the control strategy at time $t$ but only on the realization of control action at time $t.$ Thus, the future information states are {separated} from the choice of the current control strategy. This fact is crucial for deriving the optimal control strategy where, at each step, the optimization problem is to find the best control action for a given realization of the information state.
	\\
	
	The team's information structure and structural results can provide the framework to derive optimal {control prescription functions}  that will yield the optimal decisions of CAVs. Since the structural results can help us  restrict our attention to control strategies in a time-invariant domain. Thus, the optimal planning strategy of the CAVs can be derived a priori even before they start evolving in the field encountering a specific traffic scenario.
	Then, while each robotic vehicle in the team $k\in\mathcal{K}$ operates,  its control prescription function $\Gamma_t^k\colon \mathcal{L}^k\to \mathcal{U}^k$  maps the $k$'s information $\Lambda_t^k$ at time $t$ to her decision, e.g., $U_{t}^{k}=\Gamma_t^k(\Lambda_t^k)$.
	Note, the prescription functions are derived by the information structure of the team through $\Delta_t$, i.e., $\Gamma_t^{1:K}=\psi_t(\Pi_t)$.
	\\
	
	Thus, the problem for each team member $k\in\mathcal{K}$ is reformulated as to derive its {optimal planning strategy} $\psi_t^{k^*}$ so that the team's planning strategy {\boldmath$\psi_t^*$}$=$ $\{\psi_t^{k^*}; t=1,\ldots,T; k\in\mathcal{K}\}$ minimizes the expected total cost
	\begin{align}\label{problem2}
		\bar{J}({\boldmath\psi})=\mathbb{E}^{\psi}\big\{\sum_{t=1}^T c_t(\Pi_t, \Gamma_t^{1:K})\},
	\end{align}
	where $c_t(\Pi_t, \Gamma_t^{1:K})$ corresponds to travel delay of the team, and the expectation is with respect to the joint probability of the information state $\Pi_t$ and prescription functions $\Gamma_t^{1:K}$ variables designated by the choice of {\boldmath$\psi_t$}. The solution of \eqref{problem2} can be derived using standard techniques for partially observed Markov decision processes \cite{Kumar1986}. If the observation space of the system is finite, then \eqref{problem2}  has a finite-dimensional characterization. In particular, the explicit solution to \eqref{problem2}  is a piecewise linear concave function of the information state \cite{Sondik1971}.
	\\
	
	The optimal planning strategy {\boldmath$\psi_t^*$} yields the optimal control prescription function $\Gamma_t^{k^*}$ with respect to the information state $\Pi_t$, i.e., $\Gamma_t^{k^*}=\psi_t^{k^*}(\Pi_t)$. Since $U_{t}^{k}=\Gamma_t^k(\Lambda_t^k)$ and $\Gamma_t^k(\cdot)=\psi_t^{k^*}(\mathbb{P}(X_{t-n}|\Delta_t))$, the optimization problem is equivalent to the original problem \eqref{problem1a}, and hence, $g_t^k(\Delta_t, \cdot)= \Gamma_t^k(\cdot)=\psi_t^{k^*}(\Pi_t)$.
	$\Delta_t\in\mathcal{D}_t$ is the information known by all team members at time $t$, thus the prescription function $\Gamma_t^k(\cdot)$ can be derived by each robotic vehicle $k\in\mathcal{K}$ on her own {without} any {centralized} intervention. This is a key property of the proposed framework that allows each robotic vehicle $k$ to adapt its decision and thus improve the efficiency of the team in the field based on what each team member can {learn} from the information $\Lambda_t^k\in\mathcal{L}_t^k$ of her environment as discussed next. 
	\\
	
	Since we can separate the information state from the prescription function \cite{Malikopoulos2022a}, we can develop a mechanism to learn the statistic $\mathbb{P}(X_{t-n}|\Delta_t)$ using standard machine learning (ML) techniques  that will aim CAVs at adapting their decisions designated by the control prescription functions  in situations where they encounter different behaviors from what they already know about human driving. ML has been used extensively for entitling autonomy features in automotive systems \cite{Malikopoulos2011book}.   
	The optimal control prescription functions designate how each team member will coordinate and collaborate over {variable time scales} (recall that their domain is time invariant) and without any centralized (human) intervention while there exist environmental disturbances $\{W_t; t=0,\ldots,T\}$. At time $t$, all team members know the shared information $\Delta_t$, the optimal planning strategy $\psi_t^*(\cdot)$, and the optimal prescription functions $\Gamma_t^{1:K}(\cdot)$. Each CAV $k\in\mathcal{K}$ will update her optimal planing strategy according to $\hat{\psi}_t^{k^*}(\cdot)= \Psi_t^k(\Delta_t, \psi_t^{k^*}(\cdot), \Lambda_t^k)$, where $\Psi_t^k$ is a {learning function} that will use local information and observations $\Lambda_t^k$ to update the control prescription functions $\Gamma_t^{1:K}(\cdot)$. $\Psi_t^k$ will adapt the planning strategy  CAV $k$ in a way to enhance what $k$ knows about drivers' behavior. In addition,  the learning function $\Psi_t^k$ will enhance the robustness of each CAV $k$ in the presence of a contested communication environment, as described by the sequence of independent random variable $\{V_t^k: t=0,\ldots,T; k=1,\ldots,K\}$ corresponding to the noise of each team's member observation.
	The control prescription function of $k$ will then be derived according to $\hat{\Gamma}_{t}^{k^*}(\cdot)=\hat{\psi}_t^{k^*}(\mathbb{P}(X_{t-n}|\Delta_t)=\hat{\psi}_t^{k^*}(\Pi_t)$. 
	This will aim team member $k$ at feeding the realization of her local information $\Lambda_t^k$  into her updated control prescription function to derive the optimal decision, i.e., $U_{t}^{k^*}=\hat{\Gamma}_{t}^{k^*}(\Lambda_t^k)$.

	\subsection{Evaluation and Experiments in a Scaled Smart City}
	
	Ongoing research includes implementing and validating the proposed framework in the {Information and Decision Science Lab's (IDS$^3$C) (1:25)} testbed (Fig. \ref{udsscview}) \cite{chalaki2021CSM}. 
	This testbed occupies a 20 by 20 feet area and includes 50 robotic cars (both CAVs and human-driven), which can replicate real-world traffic scenarios in a small and controlled environment. 
	
	\begin{figure}
	\centering
	\includegraphics[width=0.8\linewidth, keepaspectratio]{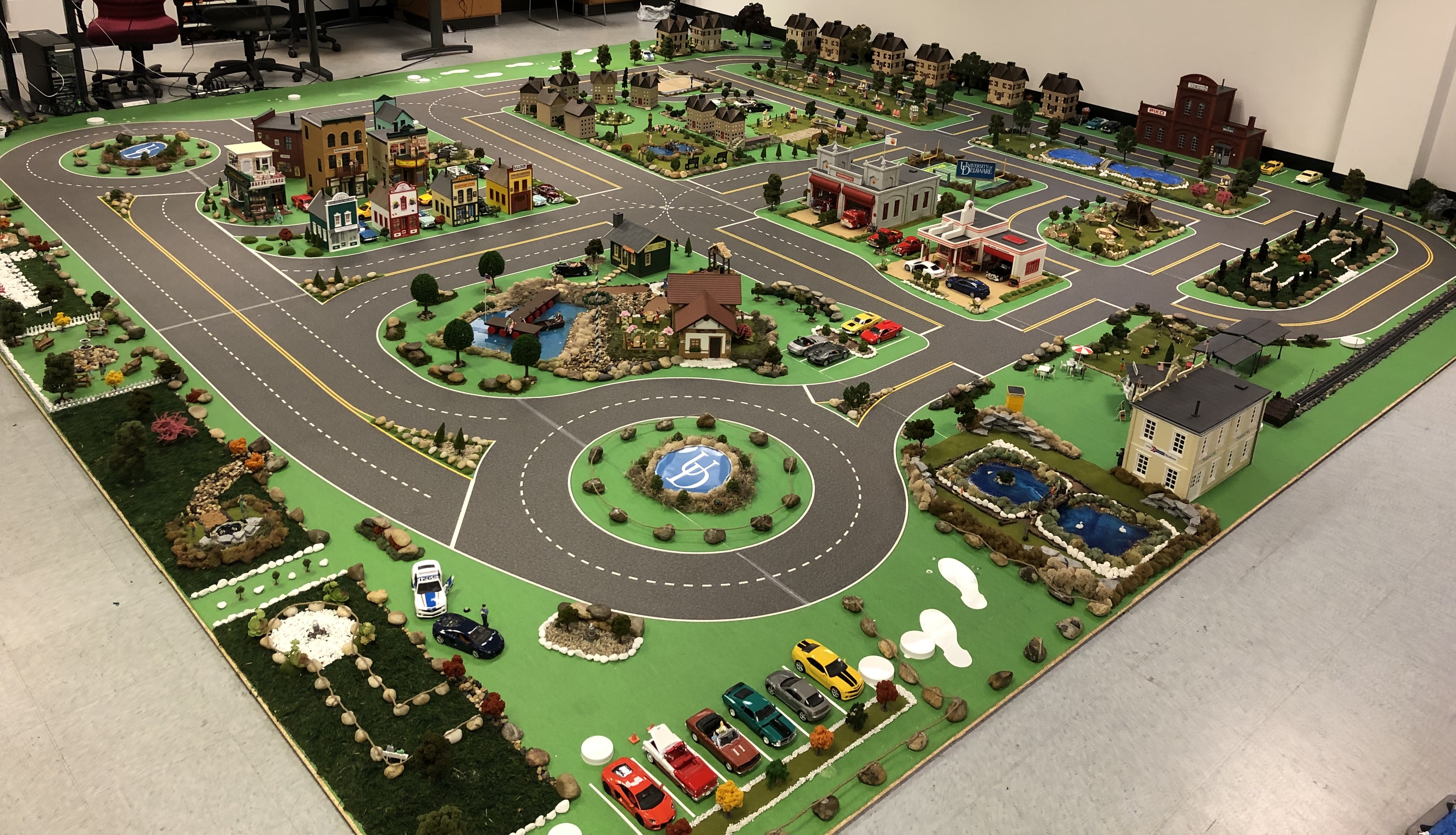} 
	\caption{A view of the IDS Lab's scaled smart city (IDS$^3$C).}%
	\label{udsscview}%
\end{figure}
	
	IDS$^3$C has 6 driver emulation stations (Fig. \ref{driver}) interfaced directly with the robotic cars which allow us to explore human driving behavior. IDS$^3$C can help us prove concepts beyond the simulation level and understand the implications of errors/delays in the vehicle-to-vehicle and vehicle-to-infrastructure communication as well as their impact on energy usage. In several recent efforts, we have used IDS$^3$C to implement and validate control algorithms for coordinating CAVs at traffic scenarios, such as merging roadways \cite{Malikopoulos2018b,Sumanth2021}, roundabouts \cite{chalaki2020experimental}, intersections  \cite{Malikopoulos2019CDC,malikopoulos2019ACC}, and corridors \cite{Beaver2020DemonstrationCity,Zhao2018ITSC,chalaki2021CSM}. We have also used IDS$^3$C to transfer policies derived in simulation using {reinforcement learning} techniques  \cite{jang2019simulation, chalaki2020ICCA}.

		\begin{figure}
		\centering
		\includegraphics[width=0.6\linewidth, keepaspectratio]{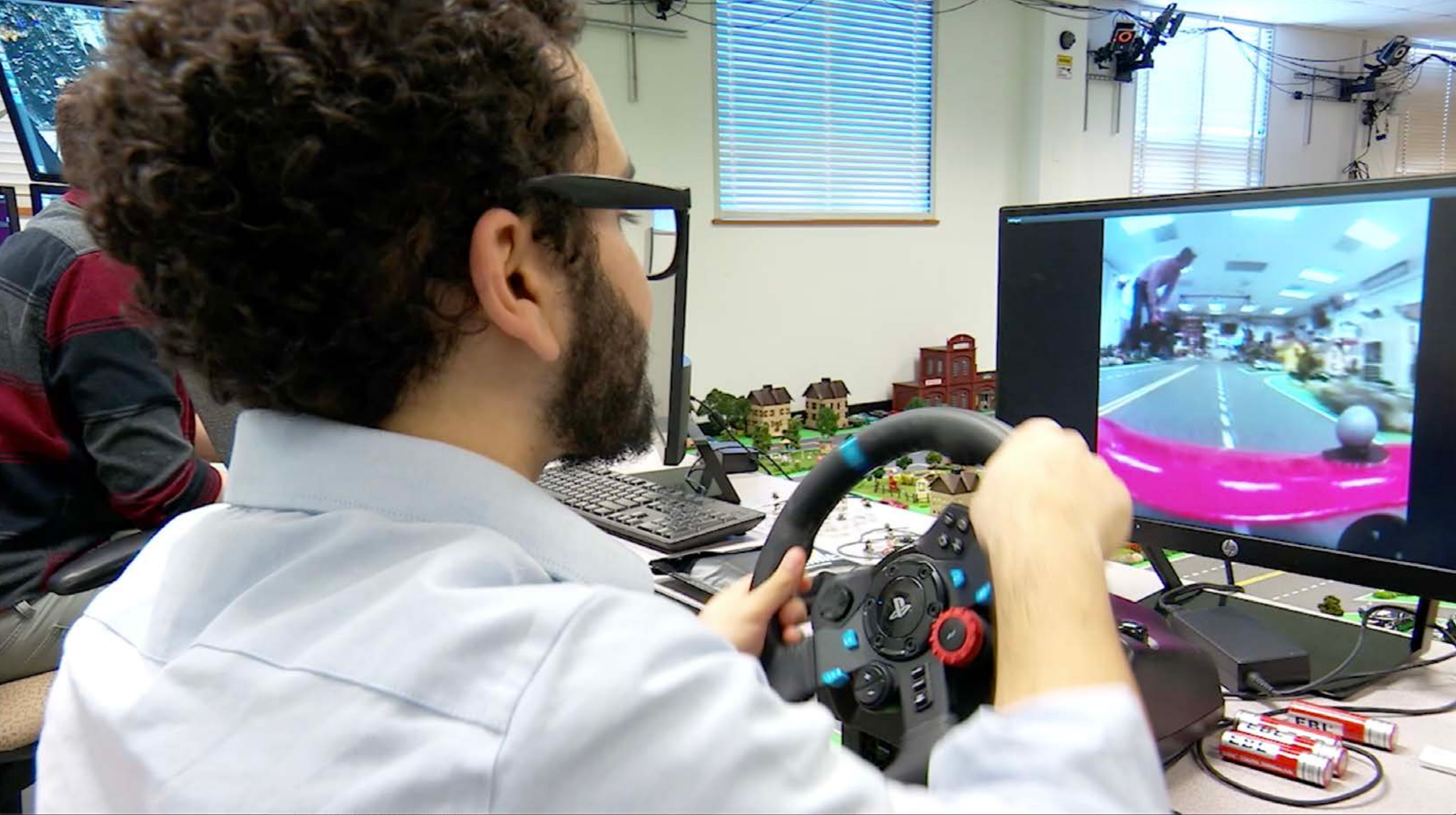} 
		\caption{A driver emulation station interfaced directly with the robotics vehicles.}%
		\label{driver}%
	\end{figure}
	
	We can apply the proposed framework  over various real-world driving scenarios deemed characteristic of typical commutes. A typical vehicle commute includes merging at roadways, crossing signalized intersections, cruising in congested traffic, and passing through speed reduction zones. The robotic cars will be able to cross intersections in  {locations 1, 3,} and {4} (Fig. \ref{udssc2}),  merge at a roundabout in  {location 2} (Fig. \ref{udssc2}), and merge at roadways in  {location 5} (Fig. \ref{udssc2}) by optimizing transportation efficiency. 
	
	\begin{figure}
	\centering
	\includegraphics[width=0.6\linewidth, keepaspectratio]{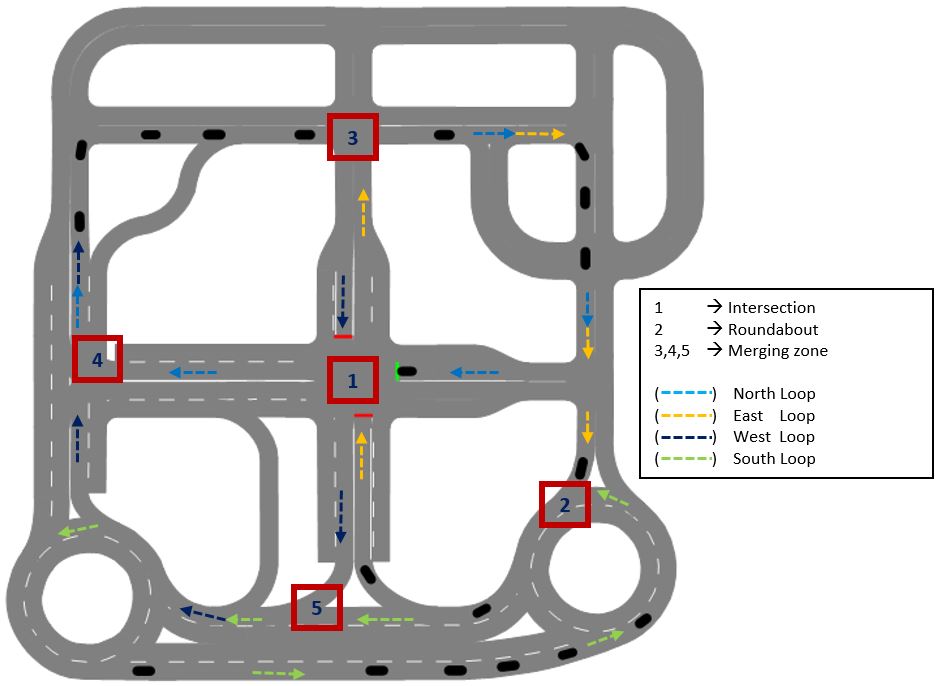} 
	\caption{Traffic scenarios of at the scaled smart city.}%
	\label{udssc2}%
\end{figure}

	\section{Concluding Remarks and Discussion}\label{sec:5}
	
	It is expected that CAVs will gradually penetrate the market and interact with HDVs in a way to improve safety and transportation efficiency \cite{Malikopoulos2018d, Zhao2018CTA, cassandras2019b} over the next several years. However, we anticipate that efficient  transportation and travel cost reduction might alter human travel behavior, causing rebound effects, e.g., by improving efficiency, travel cost is decreased, hence willingness-to-travel is increased. The latter would increase overall vehicle miles traveled, which in turn might negate the benefits in terms of energy and travel time. 
	
	We expect the proposed framework presented in this chapter to enhance our understanding of the rebound effects, travel demand and capacity changes, human reception, adoption, and use of emerging mobility systems as it addresses a complex multi-dimensional research problem focusing on societal needs. 
	The framework captures the societal impact of CAVs and provides solutions that mitigate any potential rebound effects, e.g., increased vehicle miles traveled, increased travel demand, and empty vehicle trips, while enhancing accessibility, safety, and equity in transportation.

	\section{Acknowledgments}
	This research was supported by NSF under Grants CNS-2149520 and CMMI-2219761.

\bibliographystyle{IEEEtran}
\bibliography{references_ids, tcns, references1, udssc,TAC_Ref_Andreas,MD_ref}

\end{document}